\documentclass[12pt]{article}
\usepackage{fancyhdr}
\usepackage{fullpage}
\usepackage{amsmath}
\usepackage{amssymb}
\usepackage{theorem}

\newtheorem{prop}{Proposition}[section]

{\theorembodyfont{\rm} 
 }
\newtheorem{defe}{Definition}[section]
\newtheorem{conj}{Conjecture}[section]
\numberwithin{equation}{section}
\begin{document}
\thispagestyle{plain}
\begin{center}
{\Large Verification of the GGS conjecture for $\mathfrak{sl}(n), n
\leq 12$.}

\vskip 12 pt

{\large Travis Schedler}
\end{center}
\begin{abstract}
In the 1980's, Belavin and Drinfeld classified non-unitary solutions
of the classical Yang-Baxter equation (CYBE) for simple Lie
algebras \cite{BD}.  They proved that all such solutions fall into
finitely many continuous families and introduced combinatorial objects
to label these families, Belavin-Drinfeld triples. In 1993,
Gerstenhaber, Giaquinto, and Schack attempted to quantize such
solutions for Lie algebras $\mathfrak{sl}(n).$ As a result, they
formulated a conjecture stating that certain explicitly given elements
$R \in Mat_n(\mathbb C) \otimes Mat_n(\mathbb C)$ satisfy the quantum
Yang-Baxter equation (QYBE) and the Hecke condition \cite{GGS}.
Specifically, the conjecture assigns a family of such elements $R$ to
any Belavin-Drinfeld triple of type $A_{n-1}$.  Until recently, this
conjecture has only been known to hold for $n \leq 4$.  In 1998
Giaquinto and Hodges checked the conjecture for $n=5$ by direct
computation using {\it Mathematica} \cite{GH}.  Here we report a
computation which allowed us to check that the conjecture holds for $n
\leq 12$.  The program is included which prints an element $R$ for any
triple and checks that R satisfies the QYBE and Hecke conditions.
\end{abstract}

\section{Belavin-Drinfeld triples}

Let $(e_i), 1 \leq i \leq n,$ be a basis for $\mathbb C^n$.  Set
$\Gamma = \{e_i - e_{i+1}: 1 \leq i \leq n-1\}$.  We will use the
notation $\alpha_i \equiv e_i - e_{i+1}$.  Let $( , )$ denote the
inner product on $\mathbb C^n$ having $(e_i)$ as an orthonormal basis.

\begin{defe}
A Belavin-Drinfeld triple of type $A_{n-1}$ is a triple 
$(\tau, \Gamma_1, \Gamma_2)$ where $\Gamma_1, \Gamma_2 \subset \Gamma$
and $\tau: \Gamma_1 \rightarrow \Gamma_2$ is a bijection, satisfying
two conditions:

(a) $\forall \alpha, \beta
\in \Gamma_1$, $(\tau \alpha,\tau \beta) = (\alpha, \beta)$.

(b) $\tau$ is nilpotent: $\forall \alpha \in \Gamma_1, \exists k
\in \mathbb N$ such that $\tau^k \alpha \notin \Gamma_1$.
\end{defe}

We employ three isomorphisms of Belavin-Drinfeld triples:

a) Any triple $(\tau, \Gamma_1, \Gamma_2)$ is isomorphic to the triple
$(\tau', \Gamma'_1, \Gamma'_2)$ obtained as follows: $\Gamma'_1 =
\{\alpha_m: \alpha_{n-m} \in \Gamma_1\}$, $\tau'(\alpha_m) = \alpha_k$
where $\tau(\alpha_{n-m}) = \alpha_{n-k}$.

b) Any triple $(\tau, \Gamma_1, \Gamma_2)$ is isomorphic to the triple
$(\tau^{-1}, \Gamma_2, \Gamma_1)$.

c) The product of isomorphisms (a), (b).

Modulo these isomorphisms, we found all Belavin-Drinfeld triples for
$n \leq 13$ by computer.  The number of such triples is given below:
\vskip 12 pt

\begin{center}
\begin{tabular}{|c|c||c|c||c|c|}\hline
n & \# of triples & n & \# of triples & n & \# of triples \\ \hline
2 & 1 & 6 & 41 & 10 & 10434 \\ \hline
3 & 2 & 7 & 161 & 11 & 45069 \\ \hline
4 & 4 & 8 & 611 & 12 & 201300 \\ \hline
5 & 13 & 9 & 2490 & 13 & 919479 \\ 
\hline
\end{tabular}
\end{center}

\section{The GGS conjecture}

Let $\mathfrak g = {\mathfrak sl}(n)$ be the Lie algebra of $n \times
n$ matrices of trace zero. Set $\mathfrak h \subset \mathfrak g$ to be
the subset of diagonal matrices.  Elements of $\mathbb C^n$ define
linear functions on $\mathfrak h$ by $\bigl( \sum_i \lambda_i e_i
\bigr) \bigl( \sum_i a_i e_{ii} \bigr)= \sum_i \lambda_i a_i$.  Set
$\sigma = \sum_{1 \leq i,j \leq n} e_{ij} \otimes e_{ji}$, and let $P$
be the orthogonal projection of $\sigma$ to $\mathfrak g \otimes
\mathfrak g$ with respect to the form $(X,Y) = Tr(XY)$ on
$Mat_n(\mathbb C)$.  Then, set $P^0$ to be
the projection of $P$ to $\mathfrak h \otimes \mathfrak h$.  Thus $P^0
= \sum_i \frac{n-1}{n} e_{ii} \otimes e_{ii} - \sum_{i \neq j}
\frac{1}{n} e_{ii} \otimes e_{jj}$.
 
For any Belavin-Drinfeld triple, consider the following equations:

\begin{gather} \label{r01}
r^0_{12} + r^0_{21} = P^0. \\ \label{r02} \forall \alpha \in \Gamma_1,
(\tau \alpha \otimes 1)r^0 + (1 \otimes \alpha) r^{0} = 0.
\end{gather}

Belavin and Drinfeld showed that nonunitary solutions of the CYBE
correspond to solutions of these equations.  Define $\tilde r^0 =
r^0-P^0/2$.

The GGS conjecture gives an explicit form of a matrix $R \in
Mat_n(\mathbb C) \otimes Mat_n(\mathbb C)$ for any given triple and
any given $r^0 \in \mathfrak h \otimes \mathfrak h$ satisfying
\eqref{r01}, \eqref{r02} as follows:

Set $\tilde \Gamma_1 = \{v \in \text{Span}(\Gamma_1): v = e_i - e_j, 0
\leq i < j \leq n, i \neq j\}$, and define $\tilde \Gamma_2$ similarly.
Then, extend $\tau$ to a map $\tilde \Gamma_1 \rightarrow \tilde
\Gamma_2$ so that $\tau$ is additive, i.e. $\tau(a+b) = \tau(a) +
\tau(b)$ provided $a,b,(a+b) \in \tilde \Gamma_1$.  Further, define
$\alpha \prec \beta$ if $\alpha \in \tilde \Gamma_1$ and
$\tau^k(\alpha) = \beta$, for some $k \geq 1$.  It is clear from the
conditions on $\tau$ that this means, given $\alpha = \alpha_i +
\ldots + \alpha_{i+p}$, that $\beta = \alpha_j + \ldots +
\alpha_{j+p}$, $0 \leq p \leq n-2, 1 \leq i,j \leq n, i \neq j$.
Assume $\beta = \tau^k(\alpha), k \geq 1$.  If, in this case,
$\tau^k(\alpha_i) = \alpha_{j+p}$, that is, $\tau^k$ sends the left
endpoint of $\alpha$ to the right endpoint of $\beta$, then define
$\text{sign}(\alpha,\beta) = (-1)^p$.  Otherwise, set
$\text{sign}(\alpha,\beta) = 1$.

We will use the notation $x \wedge y \equiv \frac{1}{2} (x \otimes y - y
\otimes x)$.  Furthermore, for all matrices $x \in Mat_n(\mathbb C)
\otimes Mat_n(\mathbb C)$ we will use the notation $x = \sum_{i,j,k,l}
x_{ik}^{jl} e_{ij} \otimes e_{kl}$.  Let $q$ be indeterminate and set
$\hat q \equiv q-q^{-1}$.  Finally, for any $\alpha = e_i - e_j$, set
$e_{\alpha} = e_{ij}$, and say $\alpha > 0$ if $i < j$, otherwise
$\alpha < 0$.  Now, we can define the matrix $R$ as follows:

\begin{gather} \label{ace}
a = 2 \sum_{\underset{\alpha \prec \beta}{\alpha, \beta > 0}}
\text{sign}(\alpha,\beta)\: e_{-\alpha} \wedge e_{\beta}, \quad c =
\sum_{\alpha > 0} e_{-\alpha} \wedge e_\alpha, \quad \epsilon = ac +
ca + a^2, \\ \label{tars}
\tilde a = \sum_{i,j,k,l} a_{ik}^{jl} q^{a_{ik}^{jl}
\epsilon_{ik}^{jl}}, \quad R_s = q \sum_{i} e_{ii} \otimes e_{ii} +
\sum_{i \neq j} e_{ii} \otimes e_{jj} + \hat q \sum_{i>j} e_{ij}
\otimes e_{ji}, \\ \label{r} R = q^{\tilde r^0} (R_s + \hat q \tilde a)
q^{\tilde r^0}.
\end{gather}

\begin{conj}{\bf (GGS)}  
The matrix $R$ satisfies the quantum Yang-Baxter equation, 
$R_{12} R_{13} R_{23} = R_{23} R_{13} R_{12}$, and $PR$ satisfies 
the Hecke relation, $(PR-q)(PR+q^{-1}) = 0$.
\end{conj}

\section{Checking GGS by computer}

We checked the GGS conjecture through a program written in C, which
takes as input any list of Belavin-Drinfeld triples.  For each triple,
it finds a valid $\tilde r^0$, constructs the matrix $R$, and checks
the QYBE and Hecke conditions.  Following is a more detailed
description of the procedure.

We will use the notation $\tau(\alpha) = 0$ if $\alpha \notin \tilde
\Gamma_1$.  Given a triple, the first step is to find an appropriate
$\tilde r^0$.  We rewrite the equations \eqref{r01}, \eqref{r02} as
follows:

\begin{gather} \label{tr01}
\tilde r^0_{12} + \tilde r^0_{21} = 0, \\ \label{tr02} \forall \alpha
\in \Gamma_1, ((\alpha - \tau \alpha) \otimes 1) \tilde r^0 =
\textstyle{\frac{1}{2}} ((\alpha + \tau \alpha) \otimes 1) P^0.
\end{gather}

As before, we view elements of $\mathbb C^n$ as linear functions on
$\mathfrak h$.  Then, it is easy to check $(\alpha_i)_{1 \leq i \leq
n-1}$ and $(\alpha_i - \tau \alpha_i)_{1 \leq i \leq n-1}$ are bases
of $\mathfrak h^*$.  Let $(g_i)$ and $(f_i)$ be dual to the bases
$(\alpha_i)$ and $(\alpha_i - \tau \alpha_i)$, respectively. Then, if
we view $\tilde r^0$ as an element of $Mat_{n-1}(\mathbb C)$ in the
basis $(f_i)$, it is clear that $\tilde r^0 = (b_{ij})$ where $b_{ij}
= \frac{1}{2} (\alpha_i + \tau \alpha_i, \alpha_j - \tau \alpha_j), i
\in \Gamma_1,$ where the inner product is the same we defined earlier
on $\mathbb C^n$, and $b_{ji} = -b_{ij}, i \notin \Gamma_1, j \in
\Gamma_1$.  Then, the free components of $\tilde r^0$ are those
$b_{ij}$ with $i,j \notin \Gamma_1, i < j$, which determine those
$b_{ij}, i,j \notin \Gamma_1, i > j$ since $\tilde r^0$ is
skew-symmetric.  Thus, the dimension of the space of all valid $\tilde
r^0$ is $n-m-1 \choose 2$.

The computer program merely chooses $b_{ij} = 0$ whenever $i,j \notin
\Gamma_1$.  It is known that it is sufficient to consider one element
from the family of possible $\tilde r^0$ in verifying the GGS
conjecture.  Namely, this follows from

\begin{prop} 
If $R$ of the form \eqref{r} satisfies the QYBE and PR satisfies the
Hecke relation for a given $\tilde r^0$ satisfying \eqref{tr01},
\eqref{tr02}, then for any other solution $\tilde r^0 + r'$ of
\eqref{tr01}, \eqref{tr02}, $q^{r'} R q^{r'}$ also satisfies the QYBE
and $P q^{r'} R q^{r'}$ satisfies the Hecke relation.
\end{prop}

{\it Proof.}  It is clear that $P q^{r'} R q^{r'} = q^{r'_{21}} PR
q^{r'}$. Since $r'_{21} = -r'$ by \eqref{tr01}, the Hecke relation may
be rewritten as $q^{-r'} (PR - q) (PR + q^{-1}) q^{r'} = 0$, which is
true iff $PR$ satisfies the Hecke relation.

To see that $q^{r'} R q^{r'}$ satisfies the QYBE, we take the
following steps.  By \eqref{tr02},
\begin{equation} \label{rp1}
((\alpha - \tau \alpha) \otimes 1) r' = 0.
\end{equation}
Suppose that $r' = \sum_i a_i \otimes b_i$ where the $b_i$ are
linearly independent. By \eqref{rp1}, we know that $\alpha(a_i) =
\beta(a_i)$ whenever $\alpha \prec \beta$.  Then we consider the
commutator $[a_i \otimes 1 + 1 \otimes a_i, R]$ = $[a_i \otimes 1 + 1
\otimes a_i, q^{\tilde r^0}R_s q^{\tilde r^0} + \hat q q^{\tilde r^0}
\tilde a q^{\tilde r^0}]$.  First note that $[a_i, e_{\alpha}] =
\alpha(a_i) e_\alpha$ for any $a_i \in \mathfrak h$.  Then, it is
clear $[a_i \otimes 1 + 1 \otimes a_i, q^{\tilde r^0} R_s q^{\tilde
r^0}] = [a_i \otimes 1 + 1 \otimes a_i, \sum_{i > j} d_{ij} e_{ij}
\otimes e_{ji}] = \sum_{i > j} d_{ij} (\alpha(a_i)-\alpha(a_i)) e_{ij}
\otimes e_{ji} = 0$ for the appropriate coefficients $d_{ij}$.  Now,
we see that
\begin{multline*}
[a_i \otimes 1 + 1\otimes a_i, q^{\tilde r^0} \tilde a q^{\tilde r^0}]
= [a_i \otimes 1 + 1 \otimes a_i, \sum_{\alpha, \beta > 0, \alpha
\prec \beta} (f_{\alpha,\beta} e_{-\alpha} \otimes e_{\beta} +
g_{\alpha, \beta} e_{\beta} \otimes e_{-\alpha})] \\ = \sum_{\alpha,
\beta > 0, \alpha \prec \beta} (\beta(a_i) - \alpha(a_i))
(f_{\alpha,\beta} e_{-\alpha} \otimes e_{\beta} + g_{\alpha, \beta}
e_{\beta} \otimes e_{-\alpha}) = 0.
\end{multline*}

This implies that $r' \in \Lambda^2 K$ where $K$ is the space of
symmetries of $R$, that is, $K = \{ x \in Mat_{n}(\mathbb C): [1
\otimes x + x \otimes 1, R] = 0\}$.  Furthermore, it is well-known and
easy to check that if $x \in \Lambda^2 K$ and $R$ satisfies the QYBE,
then $e^x R e^x$ also satisfies the QYBE.  Thus, in our case, we have
proved that $q^{r'} R q^{r'}$ satisfies the QYBE.  The proposition is
proved.$\quad\square$

Now, given the chosen $\tilde r^0$ in the basis $(f_i)$, the computer
program changes bases to $(g_i)$.  This is accomplished via the
transformation $[\tilde r^0]_{(g_i)} =
([(1-\tau)]_{(\alpha_i)}^{-1})^{T} [\tilde r^0]_{(f_i)}
[(1-\tau)]_{(\alpha_i)}^{-1}$ where $(1-\tau)$ is considered to be a
linear transformation on $\mathfrak h^*$, with $(1-\tau)\alpha_i =
\alpha_i - \tau \alpha_i$. Denote this new matrix by $(b'_{ij})$.

Then, the computer program obtains the matrix $[\tilde r^0]_{(e_{ii})}
\in Mat_n(\mathbb C)$ from this matrix by two quick transformations.
First it finds the intermediate matrix $(b''_{ij}) = [\tilde
r^0]_{(e_{ii}),(g_i)} \in Mat_{n \times (n-1)}(\mathbb C)$ by
$b''_{i1} = \frac{1}{n} ((n-1)b'_{i1} + (n-2)b'_{i2} + \ldots +
b'_{i,n-1})$, and the other terms follow easily.  The same technique
on the other side finally gives $[\tilde r^0]_{(e_{ii})}$.

Once $\tilde r^0$ is obtained, the computer constructs the matrix $R
\in Mat_n(\mathbb C) \otimes Mat_n(\mathbb C)$ in the basis $e_{ij}
\otimes e_{kl}, 1 \leq i,j,k,l \leq n$.  First it computes $a$, $c$,
and $\epsilon$ by \eqref{ace}.  Then, formulas \eqref{tars}, \eqref{r}
are implemented for each entry separately.  Elements $x \in
Mat_n(\mathbb C) \otimes Mat_n(\mathbb C)$ are implemented as
3-dimensional arrays $(x_{ik}^j)$, since all matrices presented in the
GGS conjecture take the form $\sum_{i,j,k} x_{ik}^j e_{ij} \otimes
e_{k,i+k-j}$.  Polynomials in $q$ are implemented as structures
containing two arrays of integers, one for positive and one for
negative powers of $q$.  The sizes of the arrays are determined in the
input of the program.

The computer checks the QYBE and Hecke conditions in the following
manner.  For the QYBE condition, the corresponding entries of $R_{12}
R_{13} R_{23}$ and $R_{23} R_{13} R_{23}$ are computed and compared
individually; both take the form
$\sum_{i,j,k,l,m} d_{ikm}^{jl} e_{ij} \otimes e_{kl} \otimes
e_{m,i+k+m-j-l}$.  The same method is applied to the Hecke condition
with matrices $\sum_{i,j,k} d_{ik}^j e_{ij} \otimes
e_{k,i+k-j}$.  Explicitly, if $R = \sum_{i,j,k} r_{ik}^j e_{ij} \otimes
e_{k,i+k-j}$, the QYBE and Hecke conditions become, respectively:

\begin{gather}
\sum_p r_{ik}^{k+i-p} r_{k+i-p,m}^j r_{p,m+k+i-p-j}^l = \sum_p
r_{km}^p r_{i,m+k-p}^{j+l-p} r_{j+l-p,p}^j, \forall i,j,k,l,m. \\
\sum_l r_{ki}^l r_{k+i-l,l}^j = \delta_{ij} + \hat q r_{ki}^j, \forall
i,j,k.
\end{gather}

Then, the computer prints the matrices $\tilde r^0$ and $R$ and
reports whether or not the conditions passed.  

After generating all Belavin-Drinfeld triples for $n \leq 13$ as
described in the previous section, all tests were performed on each
triple where $n \leq 12$ with this procedure, all of which passed.
Thus, by application of the previous proposition, we have the
following result:

\begin{prop}  
The GGS conjecture is true for Lie algebras $\mathfrak{sl}(n)$ with $n
\leq 10$.
\end{prop}

The computer program is included with this paper, with instructions on
usage included with the program itself.

\section{Acknowledgements}

I would like to thank Pavel Etingof for his generous help and advice.
I would also like to thank the Harvard College Research Program for
their support.


\begin{thebibliography}{VancVanRL}

\bibitem{BD} Belavin, A.A. and Drinfeld, V.G.: Triangle equations and
simple Lie algebras, in: S.P.Novikov (ed.), {\it Mathematical Physics
Reviews}, Harwood, New York, 1984, 93-166.

\bibitem{GGS} Gerstenhaber, M., Giaquinto, A., and Schack, S.D.:
Construction of quantum groups from Belavin-Drinfeld infinitesimals,
in: A. Joseph and S. Shnider (eds), {\it Quantum Deformations of
Algebras and their Representations,} Israel Math. Conf. Series 7,
Bar-Ilan Univ., Ramat Gan, 1993, pp. 45-64.

\bibitem{GH} Giaquinto, A. and Hodges, T.: Nonstandard solutions of
the Yang-Baxter equation, {\it Letters in Mathematical Physics} {\bf
44} (1998), 67-75.

\end{thebibliography}
\end{document}